# Obstructions to Within a Few Vertices or Edges of Acyclic

Kevin Cattell[1]   Michael J. Dinneen[1,2]   Michael R. Fellows[1]

[1] Department of Computer Science, University of Victoria,
P.O. Box 3055, Victoria, B.C. Canada V8W 3P6
[2] Computer Research and Applications, Los Alamos National Laboratory,
M.S. B265, Los Alamos, New Mexico 87545 U.S.A.

**Abstract**

Finite obstruction sets for lower ideals in the minor order are guaranteed to exist by the Graph Minor Theorem. It has been known for several years that, in principle, obstruction sets can be mechanically computed for most natural lower ideals. In this paper, we describe a general-purpose method for finding obstructions by using a bounded treewidth (or pathwidth) search. We illustrate this approach by characterizing certain families of cycle-cover graphs based on the two well-known problems: $k$-FEEDBACK VERTEX SET and $k$-FEEDBACK EDGE SET. Our search is based on a number of algorithmic strategies by which large constants can be mitigated, including a randomized strategy for obtaining proofs of minimality.

## 1  Introduction

One of the most famous results in graph theory is the characterization of planar graphs due to Kuratowski: a graph is planar if and only if it does not contain either of $K_{3,3}$ or $K_5$ as a minor. The *obstruction set* for planarity thus consists of these two graphs.

The deep results of Robertson and Seymour [8] on the well-quasi-ordering of graphs under the minor (and other) orders, have consequence of establishing non-constructively that many natural graph properties have "Kuratowski-type" characterizations, that is, they can be characterized by finite obstruction sets in an appropriate partial order. Finite forbidden substructure characterizations of graph properties have been an important part of research in graph theory since the beginning, and there are many theorems of this kind.

We describe in this paper a theory of obstruction set computation, which we believe has the potential to automate much of the theorem-proving for this kind of mathematics. This approach has successfully been used to find the obstructions for the graph families $k$-VERTEX COVER, $k = 1, \ldots, 5$ (see [1]). The pinnacle of this effort would be a computation of the obstruction set for embedding of graphs on the torus. Despite substantial human efforts for several decades this obstruction set is presently unknown (a reasonable estimate of its size seems to be about 1500 graphs).

The underlying theory for our obstruction set computations was first proved in [3], using the Graph Minor Theorem (GMT) to prove termination of the finite-state search procedure. The results in [6] can be used to prove termination without the GMT. The application of these results for the computation of any particular obstruction set requires additional problem-specific results. These results are nontrivial, but seem to be generally available (in one form or another) for virtually every natural lower ideal. Thus, this is (in principle) one route for establishing constructive versions of virtually all of the known complexity applications of the Robertson-Seymour results.

One of the curiosities of forbidden substructure theorems is the tendency of the number of obstructions for natural parameterized families of lower ideals to grow explosively as a function of the parameter $k$. For example, the number of minor-order obstructions for $k$-PATHWIDTH is 2 for $k = 1$, 110 for $k = 2$, and provably more than 60 million for $k = 3$ [5]. We favor the following as a working hypothesis: Natural forbidden substructure theorems of feasible size are feasibly computable.

The remaining sections of this paper are organized as follows. First, we formally define minor-order obstructions and the cycle-cover graph families that we characterize. Next, we present a family-independent method for finding obstruction sets. Finally, we present in the last two sections family-specific results along with obstruction sets for $k$-FEEDBACK VERTEX SET and $k$-FEEDBACK EDGE SET.

## 2 Preliminaries

Let $\leq_m$ be the *minor order* on graphs, that is, for two graphs $G$ and $H$, $H \leq_m G$ if and only if a graph isomorphic to $H$ can be obtained from $G$ by a sequence of operations chosen from: (1) taking a subgraph, and (2) contracting an edge. A family of graphs $\mathcal{F}$ is a *lower ideal* with respect to $\leq_m$ if for all graphs $G$ and $H$, the conditions (1) $H \leq_m G$ and (2) $G \in \mathcal{F}$ imply $H \in \mathcal{F}$. The *obstruction set* $\mathcal{O}_\mathcal{F}$ for $\mathcal{F}$ with respect to $\leq_m$ is the set of minimal elements of the complement of $\mathcal{F}$. This characterizes $\mathcal{F}$ in the sense that $G \in \mathcal{F}$ if and only if it is not the case that for some $H \in \mathcal{O}_\mathcal{F}$, $H \leq_m G$. The GMT states that $\leq_m$ is a well-partial order and hence $\mathcal{O}_\mathcal{F}$ is finite for all minor-order lower ideals $\mathcal{F}$.

The characterization of graph families based on the following two well-known problems (see [4]) is the focus of this paper.

**Problem 1** FEEDBACK VERTEX SET (FVS)
*Input:* Graph $G = (V, E)$ and a positive integer $k \leq |V|$.
*Question:* Is there a subset $V' \subseteq V$ with $|V'| \leq k$ such that $V'$ contains at least one vertex from every cycle in $G$?

A set $V'$ in the above problem is called a *feedback vertex set* (witness set) for the graph $G$. The family of graphs that have a feedback vertex set of size at most $k$ will be denoted by $k$-FVS. It is easy to verify that for each fixed $k$ the set of graphs in $k$-FVS is a lower ideal in the minor order. For a given graph $G$, let $FVS(G)$ denote the least $k$ such that $G$ has a feedback vertex set of cardinality $k$.

**Problem 2** FEEDBACK EDGE SET (FES)
*Input:* Graph $G = (V, E)$ and a positive integer $k \leq |E|$.
*Question:* Is there a subset $E' \subseteq E$ with $|E'| \leq k$ such that $G \setminus E'$ is acyclic?

Since the this problem closely resembles the classic FEEDBACK VERTEX SET problem, this is called the FEEDBACK EDGE SET problem. The edge set $E'$ is a *feedback edge set*. For a given graph $G$, let $FES(G)$ denote the least $k$ such that $G$ has a feedback edge set of cardinality $k$, and $k$-FES $= \{G \mid FES(G) \leq k\}$.

We conclude this section by formally defining the concept of graphs of bounded (combinatorial) width.

**Definition 3** *A* tree decomposition *of a graph $G = (V, E)$ is a tree $T$ together with a collection of subsets $T_x$ of $V$ indexed by the vertices $x$ of $T$ that satisfies:*
*1. (Covering) For every edge $(u, v)$ of $G$ there is some $x$ such that $u \in T_x$ and $v \in T_x$.*
*2. (Interpolation) If $y$ is a vertex on the unique path in $T$ from $x$ to $z$ then $T_x \cap T_z \subseteq T_y$.*
*The* width *of a tree decomposition is the maximum over the vertices $x$ of the tree $T$ of the decomposition of $|T_x| - 1$. A graph $G$ has* treewidth *at most $k$ if there is a tree decomposition of $G$ of width at most $k$.* Path decompositions *and* pathwidth *are defined by restricting the tree $T$ to be simply a path.*

# 3 How To Find Obstructions Efficiently

We search for obstructions within the set of graphs of *bounded pathwidth* (and adaptable to *bounded treewidth*). We now describe an algebraic representation for these bounded-width graphs.

**Definition 4** *A $t$-boundaried graph $G = (V, E, \partial, f)$ is an ordinary graph $G = (V, E)$ together with (1) a distinguished subset of the vertex set $\partial \subseteq V$ of cardinality $t$, the* boundary *of $G$, and (2) a bijection $f : \partial \to \{1, 2, \ldots, t\}$.*

The graphs of pathwidth at most $t$ are generated exactly by strings of (unary) operators from the following *operator set* $\Sigma_t = V_t \cup E_t$:

$$V_t = \{ \boxed{0}, \ldots, \boxed{t} \} \quad \text{and}$$
$$E_t = \{ \boxed{i\ j} \ : \ 0 \leq i < j \leq t \}.$$

To generate the graphs of treewidth at most $t$, an additional (binary) operator $\oplus$, called *circle plus*, is added to $\Sigma_t$. The semantics of these operators on $(t+1)$-boundaried graphs $G$ and $H$ are as follows:

$G\,\boxed{i}$    Add an isolated vertex to the graph $G$, and label it as the new boundary vertex $i$.

$G\,\boxed{i\ j}$    Add an edge between boundary vertices $i$ and $j$ of $G$ (ignore if operation causes a multi-edge).

$G \oplus H$    Take the disjoint union of $G$ and $H$ except that equal-labeled boundary vertices of $G$ and $H$ are identified.

A graph described by a string (tree, if $\oplus$ is used) of these operators is called a *$t$-parse*, and has an implicit labeled boundary $\partial$ of $t+1$ vertices. Throughout this paper, we refer to a $t$-parse and the graph it represents interchangeably. By convention, a $t$-parse always begins with the string $[\boxed{0}, \boxed{1}, \ldots, \boxed{t}]$ which represents the edgeless graph of order $t + 1$.

For ease of discussion throughout the remaining part of this paper, we limit ourselves to bounded pathwidth in the obstruction set search theory and only point out places where any difficulty may occur with a bounded treewidth search.

**Example 5** *A 2-parse and the graph it represents (the shaded vertices denote the final boundary).*

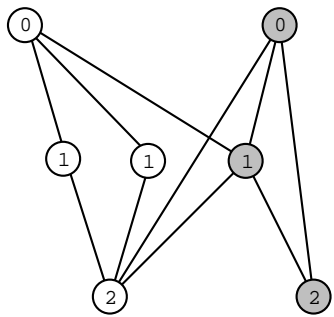

$[\boxed{0}, \boxed{1}, \boxed{2}, \boxed{0\ 1}, \boxed{1\ 2}, \boxed{1}, \boxed{0\ 1}, \boxed{1\ 2}, \boxed{1}, \boxed{0\ 1}, \boxed{1\ 2}, \boxed{0}, \boxed{0\ 1}, \boxed{0\ 2}, \boxed{2}, \boxed{0\ 2}, \boxed{1\ 2}]$

**Definition 6** Let $G = (g_1, g_2, \ldots, g_n)$ be a t-parse and $Z = (z_1, z_2, \ldots, z_m)$ be any sequence of operators over $\Sigma_t$. The concatenation $(\cdot)$ of $G$ and $Z$ is defined as

$$G \cdot Z = (g_1, g_2, \ldots, g_n, z_1, z_2, \ldots, z_m).$$

The t-parse $G \cdot Z$ is called an extended t-parse, and $Z \in \Sigma_t^*$ is called an extension. (For the treewidth case, $G$ and $Z$ are viewed as two connected subtree factors of a parse tree $G \cdot Z$ instead of two parts of a sequence of operators.)

The following sequence of definitions and results forms our theoretical basis for computing minor-order obstruction sets.

**Definition 7** Let $G$ be a t-parse. A t-parse $H$ is a $\partial$-minor of $G$, denoted $H \leq_{\partial m} G$, if $H$ is a combinatorial minor of $G$ such that no boundary vertices of $G$ are deleted by the minor operations, and the boundary vertices of $H$ are the same as the boundary vertices of $G$.

**Definition 8** Let $G$ be a t-parse. $H$ is a one-step $\partial$-minor of $G$ if $H$ is obtained from $G$ by a single minor operation (one isolated vertex deletion, one edge deletion, or one edge contraction).

Both $k$-PATHWIDTH, the family of graphs of pathwidth at most $k$, and $k$-TREEWIDTH are lower-ideals in the minor order so a $\partial$-minor $H$ of a t-parse $G$ can be represented as a t-parse. Our minor-order algorithms actually operate on the t-parses directly, bypassing any unnecessary conversion to and from the standard graph representations.

**Definition 9** Let $\mathcal{F}$ be a fixed graph family and let $G$ and $H$ be t-parses. We say $G$ and $H$ are $\mathcal{F}$-congruent (written $G \sim_{\mathcal{F}} H$) if for all extensions $Z \in \Sigma_t^*$,

$$G \cdot Z \in \mathcal{F} \iff H \cdot Z \in \mathcal{F}.$$

If $G$ is not congruent to $H$, denoted by $G \not\sim_{\mathcal{F}} H$, then we say $G$ is distinguished from $H$ (by $Z$), and $Z$ is a distinguisher for $G$ and $H$. Otherwise, $G$ and $H$ agree on $Z$.

**Definition 10** A set $T \subseteq \Sigma_t^*$ is a testset if $G \not\sim_{\mathcal{F}} H$ implies there exists $Z \in T$ that distinguishes $G$ and $H$.

In the more familiar and general setting of t-boundaried graphs (using an analogue of the Myhill-Nerode Theorem [3]), a test set $T$ may be considered to be a subset of t-boundaried graphs where concatenation $(\cdot)$ is replaced solely by circle plus $\oplus$. As we will see later, a testset is only useful for finding obstruction sets if it has finite cardinality.

**Definition 11** A t-parse $G$ is nonminimal if $G$ has a $\partial$-minor $H$ such that $G \sim_{\mathcal{F}} H$. Otherwise, we say $G$ is minimal. A t-parse $G$ is a $\partial$-obstruction if $G$ is minimal and $G \notin \mathcal{F}$.

In general, if a family $\mathcal{F}$ is a minor-order lower ideal and $G$ is $\mathcal{F}$-minimal, then for each $\partial$-minor $H$ of $G$, there exists an extension $Z$ such that

1. $G \cdot Z \notin \mathcal{F}$ and,
2. $H \cdot Z \in \mathcal{F}$.

That is, there exists a distinguisher for each possible minor $H$ of $G$.

The obstruction set $\mathcal{O}_{\mathcal{F}}$ for a family $\mathcal{F}$ is obtainable from the boundary obstruction set $\mathcal{O}_{\mathcal{F}}^{\partial}$ (set of $\partial$-obstructions) by contracting (possibly zero) edges on the boundaries of $\mathcal{O}_{\mathcal{F}}^{\partial}$, whenever the search space of width $\partial - 1$ is large enough. In our search for $\mathcal{O}_{\mathcal{F}}^{\partial}$, we must prove that each t-parse generated is minimal or nonminimal. The following two results drastically reduce the computation time required to determine these proofs.

**Lemma 12** *A $t$-parse $G$ is minimal if and only if $G$ is distinguished from each one-step $\partial$-minor of $G$. Or equivalently, $G$ is nonminimal if and only if $G$ is $\mathcal{F}$-congruent to a one-step $\partial$-minor.*

**Proof.** We proof the second statement. Let $G$ be nonminimal and suppose there exists two minors $K$ and $H$ of $G$ such that $K \leq_{\partial m} H$ and $K \sim_\mathcal{F} G$. It is sufficient to show $H \sim_\mathcal{F} G$.

For all extensions $Z \in \Sigma_t^*$, if $G \cdot Z \in \mathcal{F}$ then $H \cdot Z \in \mathcal{F}$ since $H \cdot Z \leq_{\partial m} G \cdot Z$ and $\mathcal{F}$ is a $\partial$-minor lower ideal. Now let $Z$ be any extension such that $G \cdot Z \notin \mathcal{F}$. Since $K \sim_\mathcal{F} G$, we have $K \cdot Z \notin \mathcal{F}$. And since $K \cdot Z \leq_{\partial m} H \cdot Z$, we also have $H \cdot Z \notin \mathcal{F}$. Therefore, $G$ is $\mathcal{F}$-congruent to $H$. □

**Lemma 13 (Prefix Lemma)** *If $G_n = [g_1, g_2, \ldots, g_n]$ is a minimal $t$-parse then any prefix $t$-parse $G_m$, $m < n$, is also minimal.*

**Proof.** Assume $G_n$ is nonminimal. It suffices to show that any extension of $G_n$ is nonminimal. Without loss of generality, let $H$ be a one-step $\partial$-minor of $G_n$ such that for all $Z \in \Sigma_t^*$,

$$G_n \cdot Z \in \mathcal{F} \iff H \cdot Z \in \mathcal{F}.$$

Let $g_{n+1} \in \Sigma_t$ and $G_{n+1} = G_n \cdot g_{n+1}$. Now $H' = H \cdot g_{n+1}$ is a one-step $\partial$-minor of $G_{n+1}$ such that for all $Z \in \Sigma_t^*$,

$$G_{n+1} \cdot Z = G_n \cdot (g_{n+1} \cdot Z) \in \mathcal{F} \iff H' = H \cdot (g_{n+1} \cdot Z) \in \mathcal{F}.$$

Thus, any extension of $G_n$ is nonminimal. □

The above two lemmata also hold when the circle plus operator $\oplus$ is included in $\Sigma_t$. For illustration consider the Prefix Lemma: If $G$ is a nonminimal $t$-parse with a $\mathcal{F}$-congruent minor $G'$, and $Z$ is any $t$-parse, then $(G \oplus Z)'$ is a $\mathcal{F}$-congruent minor of a nonminimal $G \oplus Z$, where we use the prime symbol to denote the corresponding minor operation done to the $G$ part of $G \oplus Z$. (The awkward notation is needed since $G' \oplus Z$ may equal $G \oplus Z$ when common boundary edges exist in both $G$ and $Z$.)

The Prefix Lemma implies that every minimal $t$-parse is obtainable by extending some minimal $t$-parse, providing a finite tree structure for the search space. In other words, the search tree may be pruned whenever a nonminimal $t$-parse is found. Since most $(t+1)$-boundaried graphs have many $t$-parse representations, we can further reduce the size of the search tree by enforcing a canonical structure on the $t$-parses considered. To do this we have to ensure that every prefix of every canonic $\partial$-obstruction (a minimal leaf of the search tree) is also canonic.

We currently use the four techniques given in Figure 1 to prove that a $t$-parse in the search tree is minimal or nonminimal. They are listed in the order that they are attempted; if one succeeds, the remainder do not need to be performed. The first three of these may not succeed, though the fourth method always will. However, if we are fortunate to have a *minimal* finite-state congruence in step 2 of Figure 1 (i.e., not a refinement of the minimum automaton for $\sim_\mathcal{F}$) then we can stop at that step since distinct final states imply the existence of an extension to distinguish the two states (and their $t$-parse representatives) of the automaton. An example of such an finite-state algorithm was used in our $k$-VERTEX COVER characterizations [1].

## 4 The FVS Obstruction Set Computation

In this section we focus on the problem-specific details for finding the $k$-FVS obstructions sets (i.e., steps 2 and 4 of Figure 1). First we describe a FVS finite-state congruence on graphs of bounded pathwidth/treewidth in $t$-parse form. Next we show how to produce complete testsets for the graph families $k$-FVS, $k \geq 0$, with respect to any boundary size $t$.

> 1. *Direct nonminimal test.* These are easily observable properties of $t$-parses that imply $t$-parses nonminimal. For any $k$-FVS family, the existence of a degree one vertex is an example of such a property.
>
> 2. *Finite-state congruence algorithm.* Such an algorithm is a refinement of the minimal finite-state (linear/tree) automaton for $\sim_\mathcal{F}$. This means that if a $t$-parse $G$ and a one-step $\partial$-minor $G'$ of $G$ have the same final state, then $G \sim_\mathcal{F} G'$, and $G$ is nonminimal. If $G$ and $G'$ have distinct final states, no conclusion can be reached.
>
> 3. *Random minor-distinguisher search.* The proof that a $t$-parse $G$ is minimal can consist of a distinguisher for each one-step $\partial$-minor $G'$ of $G$. Such distinguishers can often be easily obtained by randomly generating a sequence of operators $Z$ such that $G \cdot Z \notin F$, and then checking if $G' \cdot Z \in F$.
>
> 4. *Full test set proof.* We use a complete test set (see Definition 10) to determine if a $t$-parse $G$ is distinguished from each of its one-step $\partial$-minors. A $t$-parse $G$ is nonminimal if and only if it has a one-step $\partial$-minor $G'$ such that $G$ and $G'$ agree on every test.

Figure 1: Determining if a $t$-parse is minimal or nonminimal.

## 4.1  A FVS congruence

For a fixed $t$, let the current set of boundary vertices of a $t$-parse $G_n$ be denoted by $\partial$. Our goal is to set up a dynamic-programming congruence/automaton where the state of the $t$-parse prefix $G_{m+1}$, $m < n$, can be computed in constant time (function of $t$) from the state of the prefix $G_m$. For any subset $S$ of $\partial$, we first define the following for all prefixes $G_m$ of $G_n$.

$$F_m(S) = \left\{ \begin{array}{l} \text{least } k \text{ such that there is a FVS } V \text{ of } G_m \text{ with } V \cap \partial = S \text{ and } |V| = k \\ \text{otherwise } \infty \end{array} \right\}$$

For any witness set $V$ of $G_m$ consisting of $F_m(S)$ vertices, there is an associated *witness forest* consisting of the trees that contain at least one boundary vertex in $G_m - V$. A witness forest tells us how tight the boundary vertices are held together. Some of these forests are more concise than others for representing how vertex deletions can break up the boundary.

For two witness forests $A$ and $B$, with respect to $F_m(S)$, we say $A \leq_w B$ if the following two conditions hold:

1. For any two boundary vertices $i$ and $j$, $i$ and $j$ are connected in $A$ if and only if $i$ and $j$ are connected in $B$.

2. If for any $t$-parse extension $Z$ where there exists some non-boundary vertex $b$ of $B$ such that $(B - b) \cdot Z$ is acyclic then there exists a non-boundary vertex $a$ of $A$ such that $(A - a) \cdot Z$ is acyclic.

Also two witness forests $A$ and $B$ are equivalent if $A \leq_w B$ and $B \leq_w A$. A witness forest in reduced form (minimal number of vertices) is called a *park*.

**Lemma 14** *There are at most $3t - 3$ vertices in any park for boundary size $t$.*

**Proof.** It is easy to show that a park does not contain any non-boundary degree one vertices or any degree two vertices with non-boundary neighbors. First consider the degree two non-boundary vertices (if any). For such a vertex $v$, each of its neighbors must be a boundary vertex. After viewing $v$ and its two incident edges as a single edge between two boundary vertices, we see that at most $t - 1$ such vertices can occur. (Otherwise, a cycle would exist on the boundary.)

Now consider the remaining non-boundary vertices. Let $p$ be the number of such vertices and $e$ be the edge size of the subpark. Using the fact that the size of a forest must be strictly less than the order, we have $e < t + p$. Since the sum of the vertex degrees is twice the size, we also have $t + 3 \cdot p \leq 2 \cdot e$. Combining these inequalities while solving for $p$ we get

$$\frac{t + 3 \cdot p}{2} \leq e \leq t + p - 1, \text{ or } p \leq t - 2.$$

Summing up the boundary ($t$), the degree two vertices ($t-1$), and the degree three or more vertices ($t - 2$), shows that the order of any park can be at most $3t - 3$. □

**Corollary 15** *The total number of parks with boundary size $t$ is bounded above by $(t+1)^{t-1} \cdot 2 \cdot (2t - 1)^{2t-3}$.*

**Proof.** We transform the problem into counting the number of labeled trees. Counting separately the degree two non-boundary vertices by adding a bogus vertex to connect the trees in the park partitions of $\partial$ (assuming independent from Lemma 14), the result is then a simple application of the Cayley's Tree Formula. We have used the fact that $\sum_{i=t+1}^{2t-1} i^{i-2} \leq 2 \cdot (2t - 1)^{2t-3}$ when summing over all the possible tree orders. □

The results of the previous lemma and its corollary may be strengthened. However, these bounds are sufficient for our purposes – to show that there is a manageable (constant) number of parks (i.e., these witness sets can be used as a finite-state congruence). For each subset $S$ of the set of boundary vertices, we keep track of the parks in the following sets.

$$P_m(S) = \{P \mid P \text{ is a park of } G_m \text{ with leaves and branches over } \partial \setminus S\}$$

In the same fashion that we converted our vertex cover algorithm in [1] to a finite-state congruence for a fixed upper-bound $k$, we can use the above sets, $F_m(S)$ and $P_m(S)$ for all $S \in 2^\partial$, to construct a finite-state congruence for $k$-FVS. This is accomplished by restricting the values of $F_m(S)$ to be in $\{0, 1, \ldots, k, k+1\}$ and setting any $P_m(S) = \emptyset$ for which $F_m(S) = k + 1$; we are only interested in knowing whether or not there exists a feedback vertex set of size at most $k$. (The value of $k + 1$ acts as $\infty$.) Two $t$-parses $G_m$ and $G'_{m'}$ are congruent, $G_m \sim G'_{m'}$, if $F_m(S) = F'_{m'}(S)$ and $P_m(S) = P'_{m'}(S)$ for all $S \in 2^\partial$.

Notice that the $k$-FVS congruence $\sim$ is only a refinement of the $\mathcal{F}$-congruence $\sim_\mathcal{F}$ since $G \sim H$ implies that $G \sim_\mathcal{F} H$ but $G \not\sim H$ does not imply that $G \not\sim_\mathcal{F} H$. Thus, we will need to use a complete testset for $k$-FVS to prove $t$-parses nonminimal.

## 4.2 A complete FVS testset

Surprisingly, a finite test set for the FVS $\mathcal{F}$-congruence is easy to produce. The individual tests closely resemble the parks described above. The testset that we use consists of forests augmented with isolated triangles (and/or triangles solely attached to a single boundary vertex). Our $k$-FVS testset $T_t^k$ consists of all $t$-boundaried graphs that have the following properties:

- Each graph is a member of $k$-FVS.
- Each graph is a forest with zero or more isolated triangles, $K_3$.
- Every tree component has at least two boundary vertices.
- Every isolated triangle has at most one boundary vertex.
- Every degree one vertex is a boundary vertex.
- Every non-boundary degree two vertex is adjacent to boundary vertices.

The above restrictions on members of $T_t^k$ gives us an upper bound on the number of vertices, $|V| \leq 2t + 3(k-1)$. Hence, $T_t^k$ is a finite testset. Since this testset is based solely on $t$-boundaried graphs, it is useful for both pathwidth and treewidth searches for $k$-FVS.

**Theorem 16** *The set of $t$-boundaried graphs $T_t^k$ is a complete testset for the family $k$-FVS.*

**Proof.** Assume $G$ and $H$ are two $t$-boundaried graphs that are not $\mathcal{F}$-congruent for $k$-FVS. Without loss of generality, let $Z$ be any $t$-boundaried graph that distinguishes $G$ and $H$ with $G \oplus Z \in \mathcal{F}$ and $H \oplus Z \notin \mathcal{F}$. We show how to build a $t$-boundaried graph $T \in T_t^k$ from $Z$ that also distinguishes $G$ and $H$. Let $W$ be a set of $k$ witness vertices such that $(G \oplus Z) \setminus W$ is acyclic. From $W$, let $W_G = W \cap G$, $W_\partial = W \cap \partial$ and $W_Z = W \cap Z$. Take $T'$ to be $Z \setminus W$ plus $|W_Z|$ isolated triangles, plus $|W_\partial|$ triangles with each containing a single boundary vertex from $W_\partial$. If $T'$ contains any component $C \neq K_3$ without boundary vertices, replace it with $FVS(C)$ isolated triangles. Clearly, $G \oplus T' \in k$-FVS since $W_G$ plus one vertex from each of the non-boundary isolated triangles of $T'$ is a witness set of $k$ vertices. If $H \oplus T' \in k$-FVS then this contradicts the fact that $H \oplus Z \notin k$-FVS with $W_Z$ and $W_\partial$ and the interior witness vertices of $H$ (with respect to $H \oplus T'$). Finally, we construct a distinguisher $T \in T_t^k$ by minimizing $T'$ to satisfy the 6 properties listed above. (Note that the extension $T$ is created by not eliminating any cycles in the extension $T'$.) □

For the graph family 1-FVS on boundary size 4, the above testset consists of only 546 tests. However, for 2-FVS on boundary size 5, the above testset contains a whopping set of 14686 tests. As can be seen by the increase in the number of tests, a more compact FVS testset would be needed (if possible) before we attempt to work with boundary sizes larger than 5. The large number of tests (especially $T_5^2$) for FVS indicates why using the testset step to prove $t$-parses minimal or nonminimal is the most CPU-intensive part of our obstruction set search (and is why it is attempted last in Figure 1).

## 4.3 The $k$-FVS obstructions

We now discuss the results of our search for the 1-FVS and 2-FVS obstructions. First, we need some type of lemma that bounds the search space. The following well-known treewidth bound can be found in [7] along with other introductory information concerning the minor order and obstruction sets. We provide a proof in order to suggest how generous the bound is for the $k$-FVS obstructions, a very small subset of $(k+1)$-FVS.

**Lemma 17** *A graph in $k$-FVS has treewidth at most $k+1$.*

**Proof.** Let $G = (V, E)$ be a member of $k$-FVS and $V' \subseteq V$ be a set of $k$ witness vertices such that $G' = G \setminus V'$ is acyclic. The remaining forest $G'$ has a tree decomposition $T$ of width 1. Notice that a tree decomposition $T'$ consisting of the vertex sets of $T$ augmented as $T'_i = T_i \cup V'$ is a tree decomposition for $G$ of width $k+1$. □

**Corollary 18** *An obstruction for $k$-FVS has treewidth at most $k+2$.*

**Proof.** Let $G$ be an obstruction and $v$ any vertex of $G$. By definition, $G' = G \setminus v \in k$-FVS. Since $G'$ has a tree decomposition $T$ of width at most $k+1$, adding the vertex $v$ to each vertex of $T$ yields a tree decomposition of width at most $k+2$ for $G$. □

We now consider when the pathwidth of a $k$-FVS obstruction $G$ can be larger than the treewidth bound of $k+2$. If we attempt to build a path decomposition like the tree decompositions in the proof of Lemma 4.3, we see that for the forest $G'$ resulting by deleting an arbitrary vertex $v$ and $k$ witness vertices from $G$ has to have pathwidth at least 2. From [2] we know that the forest will contain a subdivided $K_{1,3}$ for this to happen. So, such an obstruction must have at least $1+k+7$ vertices. And for pathwidth 3, the forest has to contain one of the tree obstructions of order 22, and hence $G$ has to have at least $1+k+22$ vertices for pathwidth to be more than the treewidth plus one.

**Lemma 19** *If $O$ is an obstruction to 2-FVS and has pathwidth greater than 4, then $O$ either has at least 24 vertices or is also an obstruction to $k$-PATHWIDTH, for some $k \geq 4$.*

**Proof.** Without loss of generality, assume that the pathwidth of $O$ is 5 and is not a pathwidth obstruction. There must then exist a minor $G$ of $O$ with the same pathwidth as $O$. Since $O$ is a 2-FVS obstruction, the minor $G$ must have a feedback vertex set $V$ of cardinality 2. If the forest $G' = G \setminus V$ has pathwidth 2 or less, we can build a path decomposition of $G$ of width 4 by adding the two vertices of $V$ to the sets of a path-decomposition of $G'$ of width 2. So that leaves us with the case that $G'$ must contain a tree of pathwidth at least 3. Such a tree must have at least 22 vertices so $G$ must have at least 24 vertices. Since $O$ has the same pathwidth as $G$, the obstruction $O$ of 2-FVS must also have 24 vertices. □

Any connected obstruction to 2-FVS can not contain 3 disjoint cycles, or any degree one vertices, or any consecutive degree two vertices, so having 24 or more vertices seems unreasonable. Observe that the graph $K_5$ is an obstruction to both 2-FVS and 3-PATHWIDTH (not 4!), and that most of the $k$-PATHWIDTH obstructions have degree one vertices (and other nonminimal FVS properties), so it is unlikely that the second case of the lemma is possible. Unfortunately at this time, we have not proven the impossibility of either of these two cases. We hope, with regards to 2-FVS, that we can find a definitive proof, and avoid a treewidth 4 search.

Besides the single obstruction $K_3$ for the trivial family 0-FVS, the connected obstructions for 1-FVS and the connected obstructions for 2-FVS (pathwidth $\leq 4$) are shown in Figures 3 and 5. In our figures we have presented only the connected obstructions since any disconnected obstruction $O$ of $k$-FVS is a union of graphs from $\bigcup_{i=0}^{k-1} \mathcal{O}(i\text{-FVS})$ such that $\text{FVS}(O) = k+1$.

**Example 20** *Since $K_3$ is an obstruction for 0-FVS, and $K_4$ is an obstruction for 1-FVS, the graph $K_3 \cup K_4$ is an obstruction for 2-FVS.*

Some patterns become apparent in these two sets of obstructions such as the following easily-proven observation.

**Observation 21** *For the family $k$-FVS, the complete graph $K_{k+3}$, the augmented complete graph $A(K_{k+2})$ which has vertices $\{1, 2, \ldots, k+2\} \cup \{v_{i,j} \mid 1 \leq i < j \leq k+2\}$ and edges $\{(i,j) \mid 1 \leq i < j \leq k+2\} \cup \{(i, v_{i,j}) \text{ and } (v_{i,j}, j) \mid 1 \leq i < j \leq k+2\}$, and the augmented cycle $A(C_{2k+1})$ are obstructions.*

# 5 The FES Obstruction Set Computation

We now focus on the problem-specific details for finding the $k$-FES obstruction sets (i.e. steps 1 and 4 of Figure 1).

## 5.1 A direct nonminimal FES test

We first describe a simple graph-theoretical characterization for the graphs that are within a few edges of acyclic. This trivial result also shows that Problem 2 has a linear-time recognition algorithm.

**Theorem 22** *A graph $G = (V, E)$ with $c$ components has $FES(G) = k$ if and only if $|E| = |V| - c + k$.*

**Proof.** For $k = 0$ the result follows from the standard result for characterizing forests. If $FES(G) = k$ then deleting the $k$ witness edges produces an acyclic graph and thus $|E| = |V| - c + k$. Now consider a graph $G$ with $|V| - c + k$ edges for some $k > 0$. Since $G$ has more edges than a forest can have, there exists an edge $e$ on a cycle. Let $G' = (V, E \setminus \{e\})$. By induction $FES(G') = k - 1$. Adding the edge $e$ to a witness edge set $E'$ for $G'$ shows that $FES(G) = k$. □

Unlike $k$-FVS, it is not so obvious that the family of graphs $k$-FES is a lower ideal in the minor order. However, with the above theorem one can easily prove this.

**Corollary 23** *For each $k \geq 0$, the family of graphs $k$-FES is a lower ideal in the minor order.*

**Proof.** We show that the three basic minor operations will not increase the number of edges required to remove all cycles of a graph. An isolated vertex deletion removes both a vertex and a component at the same time, so $k$ is preserved in the formula $|E| = |V| - c + k$. For an edge deletion the number of components can increase by at most one, so with $|E|$ decreasing by one, the value of $k$ does not increase. For an edge contraction, the number of vertices decreases by one, the number of edges decrease by at least one, and the number of components stays the same, so $k$ does not increase. □

The above corollary allows us to characterize each $k$-FES family in terms of obstruction sets which we abstractly characterize below.

**Theorem 24** *A connected graph $G = (V, E)$ is an obstruction for $k$-FES if and only if $FES(G) = k+1$ and every edge contraction of $G$ removes at least two edges (i.e., the open neighborhoods of adjacent vertices overlap).*

**Proof.** This follows from the fact that an edge contraction that does not remove at least two edges is the only basic minor operation that does not decrease the number of edges required to kill all cycles for a connected graph. □

The above theorem gives us a precise means of testing for nonminimal $t$-parses (see step 1 of Figure 1).

## 5.2 A complete FES testset

Somewhat surprisingly, an usable testset for FES has already been presented in Section 4.2. We now prove that the feedback-vertex-set tests can also be used here.

**Lemma 25** *The testset $T_t^k$ for the family $k$-FVS is also a testset for $k$-FES.*

**Proof.** First observe that $k$-FES $\subseteq$ $k$-FVS so that the $k$-FVS membership restriction for $T_t^k$ graphs does not preclude any important tests (just includes some obsolete tests not in $k$-FES). Consider a fixed family $k$-FES and boundary size $t$. It suffices to show that if $G \not\sim_\mathcal{F} H$ then there exists a test $T \in T_t^k$ that distinguishes $G$ and $H$. Since $G$ and $H$ are not congruent there exists a $t$-boundaried graph $Z$ such that, without loss of generality, $G \cdot Z \in k$-FES and $H \cdot Z \notin k$-FES. We now show how to minimize $Z$ into a $T \in T_t^k$. Let $E$ be a witness edge set for $G \cdot Z \in k$-FES and let $E_Z = E(Z) \setminus E$. The first transformation on $Z$ is to set $Z' = (Z \setminus E_Z) \cup |E_Z| \cdot K_3$. Clearly $Z'$ is also a distinguisher for $G$ and $H$ since (1) $G \cdot Z' \in k$-FES by using the edges $E \setminus E_Z$ and one edge from each of the new $K_3$'s as a witness set, and (2) $H \cdot Z' \notin k$-FES, for otherwise, $H \cdot Z$ would be in $k$-FES. Notice that $Z'$ is a set of trees and isolated triangles. The final transformation on $Z$ is to let $Z''$ be $Z'$ with all non-boundary leaves deleted and non-boundary subdivided edges contracted to satisfy the conditions of a member of $T_t^k$. □

It is interesting to notice from the above proof that, in addition to the out-of-family tests, the isolated triangles in the tests for $k$-FES can be restricted to have no boundary vertices. Thus, the number of graphs in a testset for $k$-FES can be substantially smaller than our $k$-FVS testset.

## 5.3 The $k$-FES obstructions

Since the family $k$-FES $\subseteq$ $k$-FVS we know that the maximum treewidth of any obstruction for $k$-FES is at most $k + 2$. Thus, the same arguments given in Section 4.3 regarding pathwidth apply to $k$-FES as well.

For the family 0-FES, it is trivial to show that $K_3$ is the only obstructions. The connected obstructions for the graph families 1-FES through 3-FES are shown in Figures 2, 4 and 6. There are well over 100 connected obstructions for the 4-FES family. Any disconnected obstruction for $k$-FES is easily determined by combining connected obstructions from $j$-FES, $j < k$, since $FES(G_1) + FES(G_2) = FES(G_1 \cup G_2)$.

An open problem is to determine a constructive method for finding the obstructions for $k$-FES directly from the $j$-FES, $j < k$. Some easily-observed partial results are given next.

**Observation 26** *If a connected graph $G$ is an obstruction for $k$-FES then the following are all connected obstructions for $(k+1)$-FES.*

*1. $G$ with an added subdivided edge attached to an edge of $G$.*

*2. $G$ with an attached $K_3$ on one of the vertices of $G$.*

*3. $G$ with an added edge $(u,v)$ when there exists a path of length at least two between $u$ and $v$ in $G \setminus E$ for each feedback edge set $E$ of $k+1$ vertices.*

It is easy to see that if an obstruction has a vertex of degree two then it is predictable by observations (1–2). The first 2-FES obstruction in Figure 4 (wheel $W_3$) and the second 3-FES obstruction in Figure 6 ($W_4$) are two examples of graphs where observation (3) predicts the graph. Those 4-FES and 5-FES obstructions (pathwidth bound of 4) without degree two vertices and cut-vertices are shown in Figures 7 and 8. Note that the third 4-FES obstruction in Figure 7 is not predicatable from the 3-FES obstructions by using any of the above observations since deleting any edge from this graph leaves a contractable edge that does not remove any cycles (see Theorem 24).

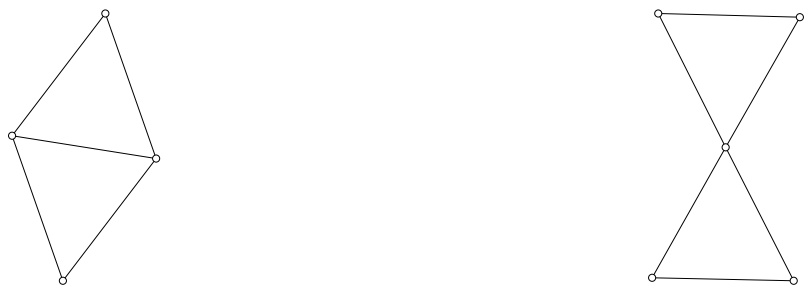

Figure 2: Connected obstructions for 1-FEEDBACK EDGE SET.

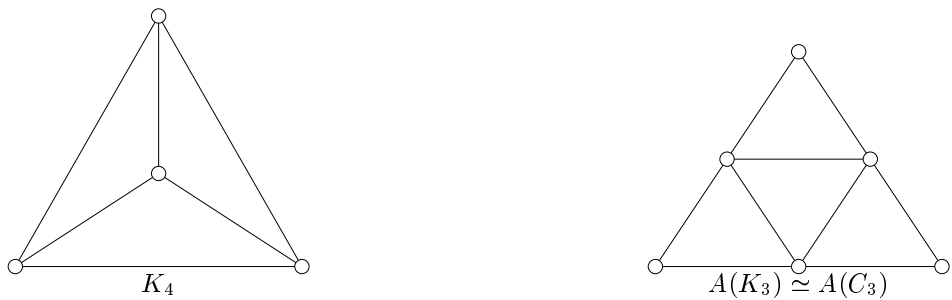

Figure 3: Connected obstructions for 1-FEEDBACK VERTEX SET.

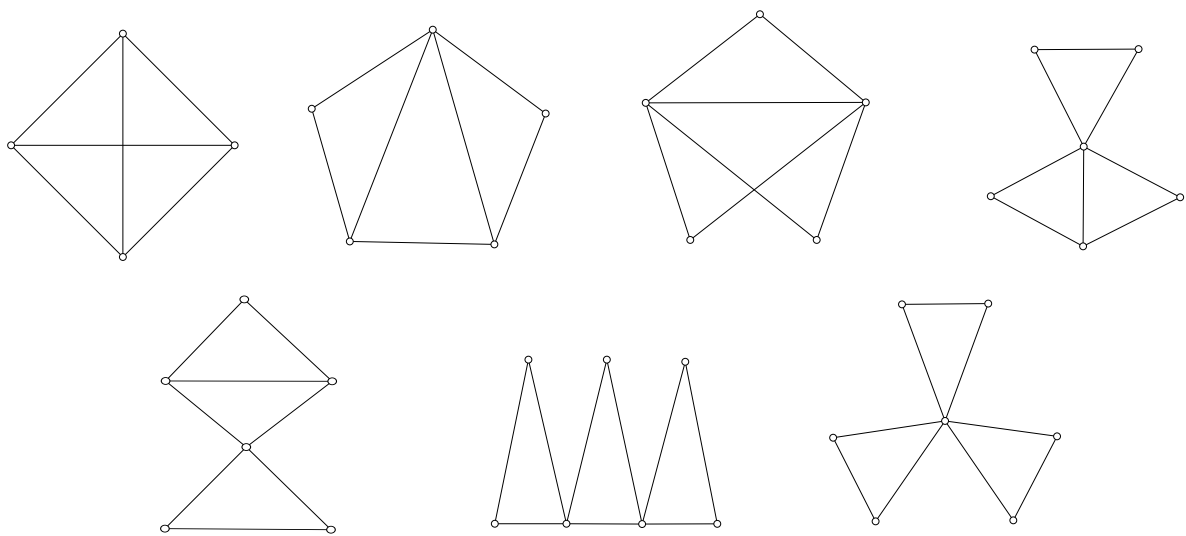

Figure 4: Connected obstructions for 2-FEEDBACK EDGE SET.

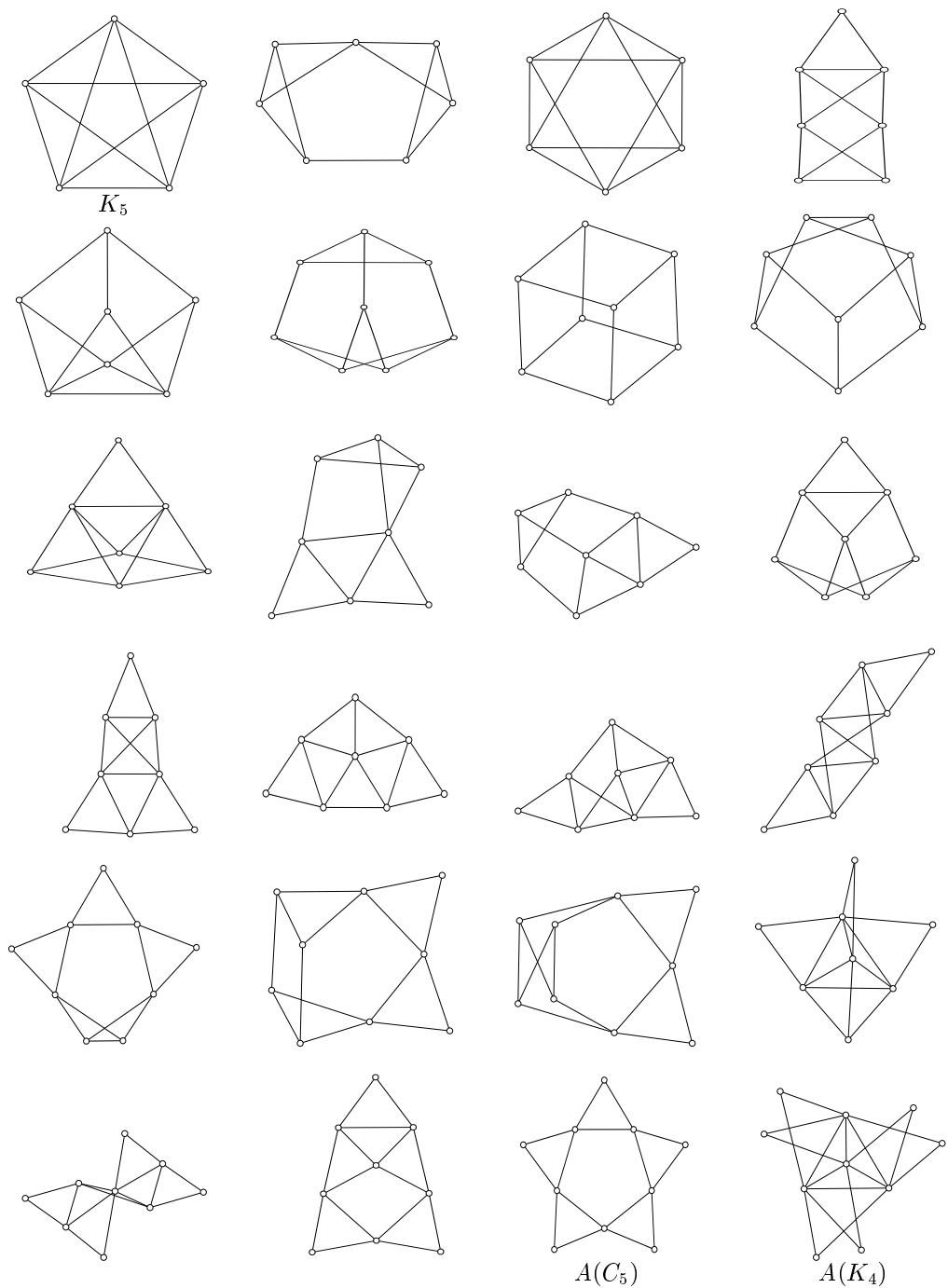

Figure 5: Connected obstructions for 2-Feedback Vertex Set, pathwidth ≤ 4.

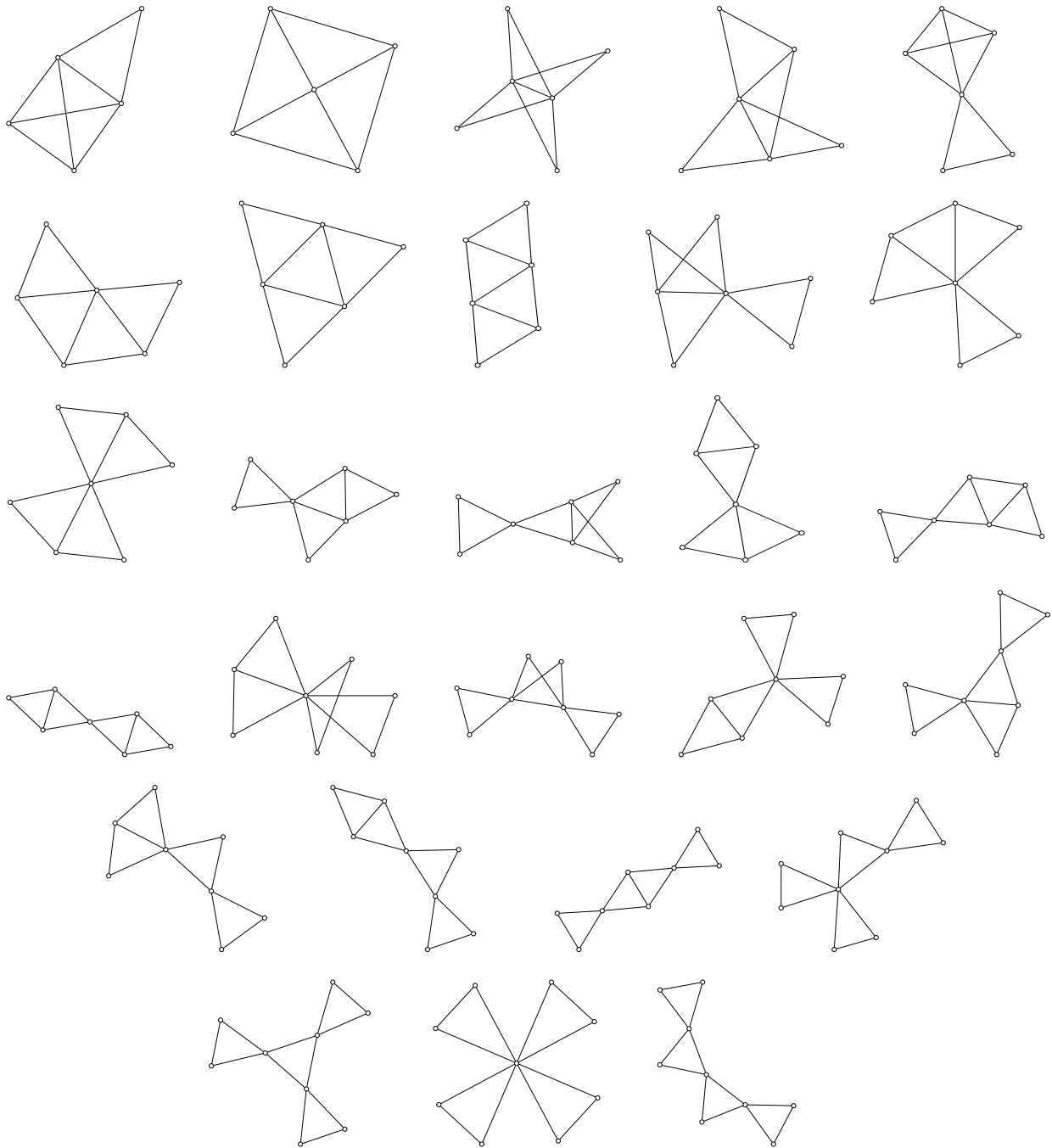

Figure 6: Connected obstructions for 3-Feedback Edge Set.

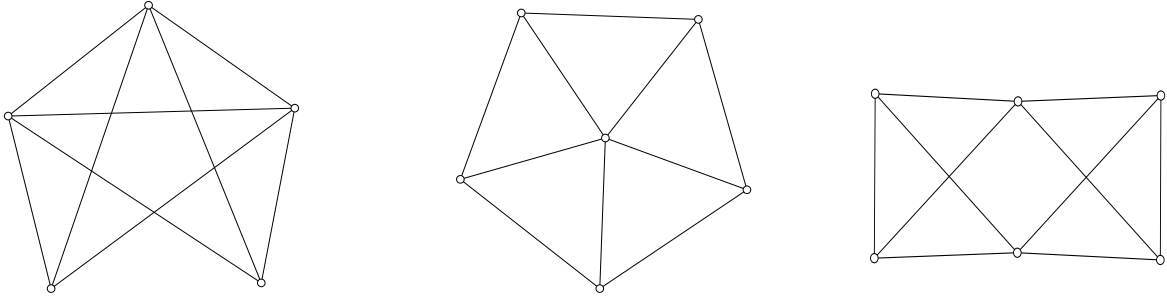

Figure 7: Biconnected 4-FES obstructions without degree 2 vertices.

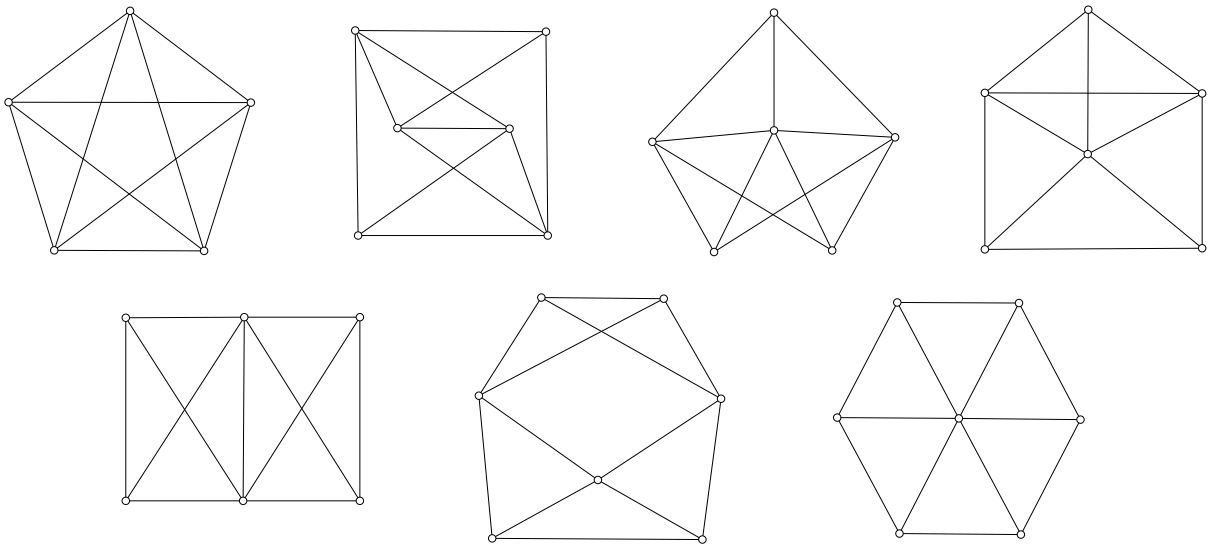

Figure 8: Biconnected 5-FES obstructions without degree 2 vertices, pathwidth $\leq 4$.